\documentclass[12pt, a4paper,twoside]{article}
\usepackage{hyperref}
\usepackage{enumerate}
\usepackage{amsmath}
\usepackage{pb-diagram}
\usepackage[english]{babel}
\usepackage{amssymb,latexsym}
\usepackage{amssymb}
\usepackage{latexsym}

\def \lra{\longrightarrow}

\newcommand{\CC}{{\mathbb C}}

\def \la{\langle}
\def \ra{\rangle}
\def \lra{\longrightarrow}
\def \lmt{\longmapsto}

\newcommand{\RR}{{\mathbb R}}
\renewcommand{\Re}{\mathrm{Re}}
\newcommand{\ZZ}{{\mathbb Z}}

\def\SU{{\rm{SU}}}
\def \U{{\rm{U}}}
\def \wnabla{\wit{\nabla}}
\def \wih{\widehat}
\def \wit{\widetilde}

\newcommand{\NN}{{\mathbb N}}

\renewcommand{\phi}{\varphi}

\newcommand{\Spin}{\mathrm{Spin}}
\renewcommand{\Re}{\mathrm{Re}}

\newcommand{\tr}{\mathrm{tr}}

\newcommand{\Spinc}{{\rm Spin}^c}

\newtheorem{example}{Examples}[section]
\newtheorem{thm}{Theorem}[section]

\newtheorem{prop}[thm]{Proposition}

\newtheorem{remark}[thm]{Remark}
\newtheorem{remarks}[thm]{Remarks}
\newtheorem{definition}[thm]{Definition}
\newtheorem{notation}[thm]{Notation}
\newtheorem{exabout:ample}[thm]{Example}

\title{{\bfseries The twisted Spin$^c$ Dirac operator on  K\"ahler submanifolds of the complex projective space}}
\author{Georges Habib\footnote{Lebanese University, Faculty of Sciences II, Department of Mathematics, P.O. Box 90656 Fanar-Matn, Lebanon,
E-mail: \texttt{ghabib@ul.edu.lb}}\, and Roger Nakad\footnote{Max Planck Institute for Mathematics, Vivatsgasse 7, 53111 Bonn, Germany, 
E-mail: \texttt{nakad@mpim-bonn.mpg.de}}}

\begin{document}
\maketitle
\begin{abstract}
\noindent In this paper, we estimate the eigenvalues of the 
twisted Dirac operator on K\"ahler submanifolds of the complex projective space $\CC P^m$ and we 
discuss the sharpness of this estimate for the embedding $\CC P^d\rightarrow \CC P^m$. 
\end{abstract} 
\section{Introduction} In his Ph.D. thesis \cite{Ginthese}, N. Ginoux gave
 an upper bound for the eigenvalues of the twisted Dirac operator for a K\"ahler spin submanifold $M^{2d}$ of a 
K\"ahler spin manifold $\wit{M}^{2m}$ carrying K\"ahlerian Killing spinors (see Equation (\ref{eq:Killing})). More precisely, he showed that 
there are  at least $\mu$ eigenvalues $\lambda_1, \lambda_2, \cdots, \lambda_\mu$ of the square of the twisted Dirac operator satisfying
\begin{eqnarray}\label{N-Gin}
 \lambda  \leqslant 
\left\{\begin{array}{ll} (d+1)^2&\textrm{ if } d \textrm{ is odd},\\\\ 
d(d+2) &\textrm{ if } d \textrm{ is even}. 
\end{array}\right.
\end{eqnarray}
Here $\mu$ denotes the dimension of the space of K\"ahlerian Killing spinors on $\wit{M}^{2m}$. Recall that 
the normal bundle is endowed with the induced spin structure coming from both manifolds $M$ and $\wit{M}$. 
The idea consists in computing the so-called Rayleigh-quotient applied to  the K\"ahlerian Killing spinor restricted to the submanifold $M$. 
The upper bound is then deduced by using the min-max principle.  This technique was also used by C. B\"ar in \cite{Bar98} for submanifolds in $\mathbb{R}^{n+1}, \mathbb{S}^{n+1}$ and $\mathbb{H}^{n+1}$.\\\\
The complex projective space $\CC P^m$ is a spin manifold if and only $m$ is odd. In this case, 
the sharpness of the upper bound (\ref{N-Gin}) was studied in \cite{GinHabibspecCPdtorduarxiv} for the canonical 
embedding $\CC P^d\rightarrow \CC P^m$, where $d$ is also odd. In fact, it is shown that for $d=1,$ the upper estimate 
is optimal for $m=3,5,7$ while it is not for $m\geq 9.$ \\\\
K\"ahler manifolds are not necessary spin but every K\"ahler manifold has a canonical  $\Spinc$ structure (see Section 2) and any  other $\Spinc$ structure can be expressed in
terms of the canonical one. Moreover,  O. Hijazi, S. Montiel and F. Urbano \cite{omu} constructed on K\"{a}hler-Einstein manifolds with positive scalar curvature, $\Spinc$ structures carrying K\"{a}hlerian
Killing spinors. Thus one can consider the result of N. Ginoux for $\Spinc$ manifolds.\\\\ 
Section 2 is devoted to recall some basic facts on $\Spinc$ structures on K\"ahler manifolds. In Section $3$, we extend the estimate (\ref{N-Gin}) to the eigenvalues of the
twisted Dirac operator for a K\"ahler submanifold of the complex projective space (see Theorem \ref{main}). 
Finally, we discuss the sharpness for the embedding $\CC P^d\rightarrow \CC P^m$ with different values of $m$ and $d$.
\section{K\"ahler Submanifolds of K\"ahler manifolds}
Let $(M^{2m}, g, J)$ be a K\"{a}hler manifold of complex dimension $m$. Recall that the complexified tangent bundle splits into the orthogonal sum $T^\CC M= T_{1,0} M \oplus T_{0,1} M$ where $T_{1,0}M$ (resp. $T_{0,1} M$) denotes
the eigenbundle of $T^\CC M$ corresponding to the eigenvalue $i$ (resp.
$-i$) of the extension of $J$. Using this decomposition, we define $\Lambda^{0, r}M:=\Lambda^r(T^*_{0,1} M)$ (resp. $\Lambda^{r,0}M$) as being the bundle of complex $r$-forms of type
$(0, r)$ (resp. of type $(r,0)$). Recall also that every K\"{a}hler manifold has a {\it canonical}
Spin$^c$ structure whose complex spinorial bundle is given  by $\Sigma M
= \Lambda^{0, *} M =\oplus_{r=0}^m \Lambda^{0, r}M,$ where the auxiliary bundle of this Spin$^c$ structure is
given by $K_M^{-1}$. Here $K_M$ is the canonical bundle of $M$ defined by $K_M= \Lambda^{m,0}M$
\cite{fr1, bm}. On the other hand, the spinor bundle of any other $\Spinc$ structure can be written as
\cite{fr1, omu}:
$$\Sigma M =  \Lambda^{0, *} M \otimes \mathfrak L,$$
 where $\mathfrak L^2 = K_M\otimes L$ and  $L$ is the auxiliary bundle
associated with this $\Spinc$ structure. Moreover, the action of the K\"ahler form $\Omega$ of $M$ splits the spinor bundle into \cite{fr1,
kirch, oussama1}:
$$\Sigma M = \oplus_{r=0}^m \Sigma_r M,$$
where $\Sigma_r M$ denotes the eigensubbundle corresponding with the
eigenvalue $i(2r-m)$ of $\Omega$ with complex rank $\Big(^m_k\Big)$. For
any vector field $X\in \Gamma(TM)$ and $\psi \in \Gamma(\Sigma_r M),$ we
have the following property $p_{\pm}(X)\cdot\psi \in \Gamma(\Sigma_{r\pm 1}M),$ where $p_{\pm}(X) = \frac 12 (X \mp iJX)$.\\\\
Let $(M^{2d},g,J)$ be an immersed K\"ahler submanifold in a K\"ahler manifold $(\wit{M}^{2m},g,J)$ carrying 
the induced complex structure $J$ (i.e. $J(TM)=TM$) and denote respectively
 by $\Omega_{\wit{M}}$, $\Omega$ and $\Omega_N$ the K\"ahler form of $\wit{M},\, M$ and of 
the normal bundle $NM\lra M$ of the immersion. Since the manifolds $M$ and ${\wit{M}}^{2n}$ are 
K\"ahler, they carry $\Spinc$ structures with corresponding auxiliary line bundles $L_M$ and $L_{\wit{M}}$. 
This induces a $\Spinc$ structure on the bundle $NM$ such that the restricted complex spinor bundle 
$\Sigma\wit{M}_{|_M}$ of $\wit{M}$ can be identified with $\Sigma M\otimes\Sigma N$, where $\Sigma M$ and 
$\Sigma N$ are the spinor bundles of $M$ and $NM$ respectively (\cite{Bar98}, \cite{GinMor02}). Moreover, 
the auxiliary line bundle $L_N$ of this $\Spinc$ structure on $NM$ is given
 by $L_N := ({L_M})^{-1} \otimes ({L_{\wit{M}}})_{|_M}.$ Given connection $1$-forms  on $L_M$ and $L_{\wit{M}}$, 
they induce a connection $\nabla:=\nabla^{\Sigma M\otimes\Sigma N}$ on $\Sigma M\otimes\Sigma N$. Thus one 
can state a Gauss-type formula for the spinorial Levi-Civita connections $\wnabla$ and $\nabla$ on $\Sigma\wit{M}$ and 
$\Sigma M\otimes\Sigma N$ respectively \cite{r2}. That is, for all $X\in TM$ and $\varphi\in\Gamma(\Sigma\wit{M}_{|_M})$, 
we have
\begin{equation}\label{eq:Gauss}\wnabla_X\varphi=\nabla_X\varphi+\frac{1}{2}\sum_{j=1}^{2d}e_j\cdot II(X,e_j)\cdot\varphi, \end{equation}
where $(e_j)_{1\leq j\leq 2d}$ is any 
local orthonormal basis of $TM$ and $II$ is the second fundamental form of the immersion. As a consequence 
of the Gauss formula, the square of the auxiliary Dirac-type operator $\wih{D}:=\sum_{j=1}^{2d}e_j\cdot\wnabla_{e_j}$ 
is related to the square of the twisted Dirac operator $D_M^{\Sigma N}:=\sum_{j=1}^{2d}e_j\cdot\nabla_{e_j}$ by \cite[Lemme 4.1]{Ginthese}:
$$\wih{D}^2 \phi = (D_M^{\Sigma N})^2\phi -d^2\vert H\vert^2 \phi -d\sum_{j=1}^{2d} e_j\cdot\nabla^N_{e_j}H\cdot\phi,$$
where $H := \frac{1}{2d} \tr(II)$ is the mean curvature vector field of the immersion. In our case, the mean curvature vanishes which means that the operators $\wih{D}^2$ and $(D_M^{\Sigma N})^2$ coincide.\\\\
In the sequel, take the manifold $\wit{M}$ as the complex projective space $\CC P^m$ endowed with its Fubini-Study metric of constant holomorphic sectional curvature $4$. In \cite{omu}, the authors proved that for every $q\in \ZZ$, such that $q +m+1 \in 2\ZZ$, there exists a $\Spinc$ structure on $\CC P^m$ whose auxiliary line
 bundle is given by 
$\mathcal{L}_m^q$. Here $\mathcal{L}_m$ denotes the tautological bundle of $\CC P^m.$ 
In particular for $q=-m-1$ (resp. $q=m+1$), the $\Spinc$ structure is the canonical one (resp. anti-canonical) \cite{A1} and 
for $q=0$ it corresponds to the unique spin structure if $m$ is odd. Let us denote by $\Sigma^q {\CC P^m}$ the spinor bundle of the corresponding $\Spinc$ structure with $\mathcal{L}^q$ as auxiliary line
bundle. For any integer  $r$ in $\{0,\cdots,m+1\}$ such that  $q= 2r-(m+1)$, the bundle $\Sigma^q {\CC P^m}$ carries a K\"{a}hlerian Killing spinor field $\psi = \psi_{r-1} + \psi_r$ 
satisfying, for all $X\in \Gamma(T {\CC P^m})$ \cite{omu}
\begin{eqnarray}\label{eq:Killing}
\wnabla_X \psi_r &=& - p_+(X)\cdot\psi_{r-1},\nonumber\\
\wnabla_X \psi_{r-1} &=& - p_-(X)\cdot\psi_{r},
\end{eqnarray}
The space of K\"{a}hlerian Killing spinors is of rank $\left(\begin{array}{c}m+1\\r\end{array}\right)$. We point out that for $r=0$ (resp. $r=m+1$) the K\"{a}hlerian Killing spinor is a parallel spinor which is carried by the canonical structure (resp. anti-canonical). Moreover, for $r=\frac{m+1}{2}$, i.e. $m$ is odd, the K\"{a}hlerian Killing spinor is the usual one lying in $\Sigma^0_{\frac{m-1}{2}} \CC P^m\oplus \Sigma^0_{\frac{m+1}{2}}\CC P^m$ defined in \cite{oussama1, kirch}. 
\section{Main result}
In this section, we will establish the estimates for the eigenvalues of the twisted Dirac operator of complex submanifolds of the complex projective space. We will test the sharpness of Inequality (\ref{eq:estimate}) for the canonical embedding 
$\CC P^d\rightarrow \CC P^m$. For more details, we refer to \cite{GinHabibspecCPdtorduarxiv}.
\begin{thm}\label{main}
Let $(M^{2d},g,J)$ be a closed K\"ahler submanifold of the complex projective space $\CC P^m$. For 
$r\in \{0,\cdots,m+1\}$ and $q=2r-(m+1)$, there are at least $\left(\begin{array}{c}m+1\\r\end{array}\right)$-eigenvalues $\lambda$ of $(D_M^{\Sigma N})^2$ satisfying 
\begin{equation}\label{eq:estimate}
\lambda \leqslant 
\left\{\begin{array}{ll} -(q^2-(d+1)^2) +2|q|(m-d)-1 &\textrm{ if } m-d \textrm{ is odd}\\\\ 
-(q^2-(d+1)^2) +2|q|(m-d) &\textrm{ if } m-d \textrm{ is even}. 
\end{array}\right.
\end{equation}
\end{thm}
{\bf Proof.} The proof relies on computing the Rayleigh-quotient 
$$\frac{\int_M \Re \la(D_M^{\Sigma N})^2\psi,\psi\ra v_g}{\int_M |\psi|^2 v_g}$$ applied to any non-zero K\"ahlerian Killing spinor $\psi=\psi_{r-1}+\psi_r$ on $\CC P^m$. 
A straightforward computation of the auxiliary Dirac operator leads to  
\[\wih{D}\psi_{r-1}=(q+d+1)\psi_r+i\Omega_N\cdot\psi_r.\]
\[\wih{D}\psi_{r}=-(q-d-1)\psi_{r-1}-i\Omega_N\cdot\psi_{r-1}.\]
Summing up the above two equations, we deduce after using the fact that the auxiliary Dirac operator commutes with the normal K\"ahler form \cite{GinHabibspecCPdtorduarxiv}, that 
$$\wih{D}^2 \psi = - (q^2 - (d+1)^2) \psi - 2iq \Omega^N\cdot\psi + \Omega^N\cdot\Omega^N\cdot\psi.$$
Taking the Hermitian inner product with $\psi$ and using the fact that the seond term can be bounded from above by $2|q|(m-d)$, we get our estimates after using $\vert\Omega^N\cdot\psi\vert\geq \vert \psi\vert$ if $m-d$ is odd and $0$ otherwise. 
\hfill$\square$\\\\

\noindent In the following, we will treat the sharpness through the embedding $\CC P^d\rightarrow \CC P^m$. Recall first that the complex projective 
space $\CC P^d$ can be seen as the symmetric space $\SU_{d+1}/{\rm S}(\U_d\times\U_1)$ where 
${\rm S}(\U_d\times\U_1):=\{\left(\begin{array}{ll}B&0\\0&\det(B)^{-1}\end{array}\right)\,|\,B\in\U_d\}$. The tangent 
bundle of $\CC P^d$ can be described as a homogeneous bundle which is associated with the
 ${\rm S}(\U_d\times\U_1)$-principal bundle $\SU_{d+1}\rightarrow \CC P^d$ via the isotropy representation  
\begin{eqnarray*} 
\alpha:\ \ \ \ \ \ \ \ {\rm S}(\U_d\times\U_1)&\lra&\U_d\\
\left(\begin{array}{cc}B&0\\0&\det(B)^{-1}\end{array}\right)&\lmt&\det(B) B.
\end{eqnarray*}
For the canonical embedding $\CC P^d\rightarrow \CC P^m$, the normal bundle $T^\perp{\CC P^d}$ is isomorphic 
to $\mathcal{L}_d^*\otimes \mathbb{C}^{m-d}$ where $\mathcal{L}_d$ is the tautological bundle of $\CC P^d$. 
The bundle $\mathcal{L}_d$ is isomorphic to the homogeneous bundle which is associated with the 
${\rm S}(\U_d\times\U_1)$-principal
 bundle $\SU_{d+1}$ via the representation  
\begin{eqnarray*} 
\rho:\ \ \ \ \ \ \ \ \ \ {\rm S}(\U_d\times\U_1) &\lra& \U_{1}\\
\left(\begin{array}{cc}B& 0\\0&\det(B)^{-1}\end{array}\right)&\lmt&(\det(B))^{-1}.
\end{eqnarray*}
Thus the normal bundle is associated with the ${\rm S}(\U_d\times\U_1)$-principal bundle $\SU_{d+1}\rightarrow \CC P^d$ via the representation
\begin{eqnarray*}
\rho:\ \ \ \ \ \ \ \ \ \ {\rm S}(\U_d\times\U_1)&\lra&\U_{m-d} \\
\left(\begin{array}{cc}B&0 \nonumber\\0&\det(B)^{-1}\end{array}\right)&\lmt&\det(B)\mathrm{I}_{m-d}. \nonumber
\end{eqnarray*} 
Now, we endow $\CC P^d$ with a  $\Spinc$  structure whose auxiliary line bundle is given by $\mathcal L_d^{q'}$ for $q'\in \mathbb{Z}$. In this case, its spinor bundle is given by 
$$\Sigma^{q'}\CC P^d = \Lambda^{0, *} \CC P^d\otimes \mathcal{L}_d^{\frac{q'+d+1}{2}}.$$ 
The existence of $\Spinc$ structures on both $\CC P^d$ and $\CC P^m$ induces also a $\Spinc$ structure on the normal bundle of the embedding with auxiliary line bundle is given by 
$\mathcal{L}_m^q|_{\CC P^d} \otimes \mathcal{L}_d^{q'}$ which is isomorphic to  $\mathcal{L}_d^{q-q'}$. Therefore the Lie-group homomorphism
\begin{eqnarray*} 
\rho:\ \ \ \ \ \ \ \ \ \  {\rm S}(\U_d\times\U_1) &\lra&\U_{m-d}\times {\rm U}_1\\
\left(\begin{array}{cc}B&0\\0&\det(B)^{-1}\end{array}\right)&\lmt& (\det(B)\mathrm{I}_{m-d},\det(B)^{-(q-q')})
\end{eqnarray*}
can be lifted through the non-trivial two-fold covering map $\Spinc_{2(m-d)}\lra {\rm SO}_{2(m-d)}\times {\rm U}_1$ to the homomorphism  
\begin{eqnarray*}
\wit{\rho}:\ \ \ \ \ \ \ \ \ \ {\rm S}(\U_d\times\U_1)&\lra&\Spinc_{2(m-d)}\\
\left(\begin{array}{cc}B&0\\0&\det(B)^{-1}\end{array}\right)&\lmt&(\det(B))^{-\frac{q - q'+ m-d}{2}}j(\det(B){\rm I}_{m-d}),\nonumber
\end{eqnarray*} 
where for any positive integer $k$, we recall that $j:{\rm U}_k\lra \Spinc_{2k}$ is given on elements
 of diagonal form of ${\rm U}_k$  as 
$$j({\rm diag}(e^{i\lambda_1},\cdots,e^{i\lambda_k}))=e^{\frac{i}{2}(\sum_{j=1}^k\lambda_j)}\widetilde{R}_{e_1,Je_1}(\frac{\lambda_1}{2})\cdots  \widetilde{R}_{e_k,Je_k}(\frac{\lambda_k}{2}).$$ 
Here $J$ is the canonical complex
 structure on $\CC^k$ and $\widetilde{R}_{v,w}(\lambda)={\rm cos}(\lambda)+{\rm sin}(\lambda)v\cdot w\in \Spin_{2k}$ 
is defined for any orthonormal system $\{v,w\}\in \RR^{2k}.$ We point out that the integer $q-q'+m-d$ is even. 
Following the similar proof as in \cite[Corollary 4.4]{GinHabibspecCPdtorduarxiv}, the complex spinor bundle of $T^\perp \CC P^d$ splits into the orthogonal sum
$$\Sigma(T^\perp\CC P^d)\cong\bigoplus_{s=0}^{m-d}\left(\begin{array}{c}m-d\\ s\end{array}\right)\mathcal{L}_d^{\frac{q-q'+m-d}{2}-s},$$
where for each $s\in\{0,\ldots,m-d\}$, the factor $\left(\begin{array}{c}m-d\\ s\end{array}\right)$ stands 
the multiplicity which the line bundle $\mathcal{L}_d^{\frac{q-q'+m-d}{2}-s}$ appears in the splitting. This gives the following decomposition 
\begin{eqnarray*}
\Sigma \CC P^d\otimes \Sigma^\perp \CC P^d&\simeq & \mathop\oplus\limits_{s=0}^{m-d}\begin{pmatrix}m-d\\s\end{pmatrix}\Sigma \CC P^d\otimes \mathcal{L}_d^{\frac{q-q'+m-d}{2}-s}\\
&\simeq &\mathop\oplus\limits_{s=0}^{m-d}\begin{pmatrix}m-d\\s\end{pmatrix}\Lambda^{0,*}\CC P^d\otimes \mathcal{L}_d^{\frac{d+1+q'}{2}}\otimes \mathcal{L}_d^{\frac{q-q'+m-d}{2}-s}\\
&\simeq &\mathop\oplus\limits_{s=0}^{m-d}\begin{pmatrix}m-d\\s\end{pmatrix}\Lambda^{0,*}\CC P^d\otimes\mathcal{L}_d^{\frac{m+1+q}{2}-s}.
\end{eqnarray*}  
We point out here that the above decomposition does not depend on the $\Spinc$ structure chosen on $\CC P^d$, since no power in $q'$ appears. In \cite{Murray}, the authors proved that (see also \cite{GinHabibspecCPdtorduarxiv, BenHalima08}): 
\begin{prop}\label{spe}  
Let $\CC P^d$ be the complex projective space of constant holomorphic sectionnal curvature $4$ endowed with a $\Spinc$ structure 
whose spinor bundle is given by $\Lambda^{0,*} \CC P^d\otimes \mathcal{L}^v$, for some $v\in \ZZ$, i.e., whose auxiliary line bundle is given by 
$\mathcal{L}_d^{2v-(d+1)}$. Then, the spectrum of the square of the Dirac operator is given by the eigenvalue $0$ if 
$v \leq 0$ or $v \geq d+1$ and by 
$$\lambda^2 = 4(l+v)(l-k+d),$$
where $l\in \NN$, $l+v \geq k+1$ and $0 \leqslant k \leqslant d-1$. Moreover, the multiplicity of $\lambda^2$ is given by 
$$\frac{2(l+d)!(l+v-k-1+d)! (2l+v-k+d)}{l!k!d!(l+v-k-1)!(d-k-1)! (l+v)(l+d-k)}$$
and the multiplicity of $0$ by $\frac{(\vert v\vert +d)!}{d!\vert v \vert !}$ if $v\leqslant 0$ and by 
$\frac{(v-1)!}{d!(v-d-1)!}$ if $v \geq d+1$. 
\end{prop}
In order to find the spectrum of the square of the twisted Dirac operator corresponding with the embedding  $\CC P^d\rightarrow \CC P^m$, 
one should replace $v$ in Proposition \ref{spe} by $\frac{m+1+q}{2}-s$ and  in this case, 
the eigenvalue of the square of the twisted Dirac operator is given by 
 $0$  if $\frac{m+1+q}{2}-s \leqslant 0$ or  if $\frac{m+1+q}{2}-s \geq d+1$ and by
$4 (l+\frac{m+1+q}{2} -s)(l-k+d)$ for $0 \leqslant s \leqslant m-d$, $0 \leqslant k \leqslant  d-1$ and
 $l+ \frac{m+q+1}{2} -s \geq k+1$. \\\\
Let us consider particular values for $d, m$ and $q$ in order to check the optimality. For $d=1$, $m=2$ and $q=1$, by Theorem \ref{main},  there are at least $3$ eigenvalues of the square of the twisted
Dirac operator satisfying the estimate $\lambda \leqslant 4$. The multiplicity of zero is $1$ and the multiplicity 
of the eigenvalue  $4$ is $4$ which means that the estimate is optimal. For $d=1$, $m=3$ and $q=2$, by Theorem \ref{main},  
there are at least $4$ eigenvalues of the square of the twisted
Dirac operator satisfying the estimate $\lambda \leqslant 8$. The multiplicity of zero is $3$. The multiplicity 
of the eigenvalue  $4$ is $4$ and the multiplicity of the eigenvalue $8$ is $6$ which means that the estimate is not optimal. 
For $d=2$, $m=3$ and $q=2$,   
there are at least $4$ eigenvalues of the square of the twisted
Dirac operator satisfying the estimate $\lambda \leqslant 8$. The multiplicity of zero is $1$ and the multiplicity 
of the eigenvalue $8$ is $6$ which means that the estimate is optimal.\\\\
{\bf Acknowledgment}\\\\ 
The authors would like to thank Nicolas Ginoux and Oussama Hijazi for fruitful discussions during the preparation of this work. They also thank the Max Planck Institute for Mathematics and the University of Nancy for their support.

\end{document}